\documentclass{amsart}
\usepackage{graphicx}
\usepackage{amssymb}
\usepackage{amsfonts}
\swapnumbers
\newtheorem{theorem}{Theorem}[section]
\newtheorem{lemma}[theorem]{Lemma}

\theoremstyle{definition}

\newtheorem{assumption}[theorem]{Assumption}
\newtheorem{remark}[theorem]{Remark}

\numberwithin{equation}{section}
 \theoremstyle{plain}
 
 \numberwithin{equation}{section} 
 \numberwithin{figure}{section} 
 \theoremstyle{plain}
 \theoremstyle{plain}
 \theoremstyle{remark}
 \newtheorem*{acknowledgement*}{Acknowledgement}

\newcommand{\cB}{{\mathcal B}}

\newcommand{\cF}{{\mathcal F}}

\newcommand{\cL}{{\mathcal L}}

\newcommand{\cX}{{\mathcal X}}
\newcommand{\cY}{{\mathcal Y}}

\newcommand{\te}{{\theta}}

\newcommand{\Om}{{\Omega}}
\newcommand{\om}{{\omega}}
\newcommand{\ve}{{\varepsilon}}
\newcommand{\del}{{\delta}}

\newcommand{\gam}{{\gamma}}
\newcommand{\Gam}{{\Gamma}}
\newcommand{\vf}{{\varphi}}

\newcommand{\sig}{{\sigma}}
\newcommand{\al}{{\alpha}}
\newcommand{\be}{{\beta}}
\newcommand{\ka}{{\kappa}}
\newcommand{\la}{{\lambda}}

\newcommand{\bbR}{{\mathbb R}}

\newcommand{\bbZ}{{\mathbb Z}}
\newcommand{\bbI}{{\mathbb I}}

\begin{document}
\title[]{A Nonconventional Local Limit Theorem}%
 \vskip 0.1cm
 \author{Yeor Hafouta and Yuri Kifer\\
\vskip 0.1cm
 Institute  of Mathematics\\
Hebrew University\\
Jerusalem, Israel}%
\address{
Institute of Mathematics, The Hebrew University, Jerusalem 91904, Israel}
\email{yeor.hafouta@mail.huji.ac.il, kifer@math.huji.ac.il}%

\thanks{ }
\subjclass[2000]{Primary: 60F05 Secondary: 60J05}%
\keywords{local limit theorem, Markov chain, mixing, nonconventional setup.}%
\dedicatory{  }
 \date{\today}
\begin{abstract}\noindent
Local limit theorems have their origin in the classical De Moivre--Laplace
theorem and they study the asymptotic behavior as $N\to\infty$ of probabilities
having the form $P\{ S_N=k\}$ where $S_N=\sum^N_{n=1}F(\xi_n)$ is a sum of an
integer valued function $F$ taken on i.i.d. or Markov dependent sequence of
random variables $\{\xi_j\}$. Corresponding results for lattice valued and
general functions $F$ were obtained, as well. We extend here this type of
results to nonconventional sums of the form $S_N=\sum_{n=1}^NF(\xi_n,\xi_{2n},
...,\xi_{\ell n})$ which continues the recent line of research studying
various limit theorems for such expressions.
\end{abstract}
\maketitle
\markboth{Y. Hafouta and Y. Kifer }{Local Limit Theorem }
\renewcommand{\theequation}{\arabic{section}.\arabic{equation}}
\pagenumbering{arabic}

\section{Introduction}\label{sec1}\setcounter{equation}{0}

The classical De Moivre-Laplace theorem says that if $\xi_1,\xi_2,...$
are independent identically distributed (i.i.d.) 0--1 Bernoulli random
variables taking on 1 with probability $p$ and $S_N=\sum^N_{n=1}\xi_n$
then the probability $P\{ S_N=k\}$ is equivalent as $N\to\infty$ to
$(2\pi Npq)^{-1/2}\exp(-(k-Np)^2/2Npq)$, $q=1-p$ uniformly in $k$ such
that $|k-Np|=o(Npq)^{2/3}$. It turns out that this type of results can be
extended to "nonconventional" sums of the form
\[
S_N=\sum_{n=1}^N\xi_n\xi_{2n}\cdots\xi_{\ell n}
\]
where $\{\xi_j\}$ are either an i.i.d. random variables or they form a
sufficiently fast mixing Markov chain.

In fact, we will deal here with more general sums of the form
\begin{equation}\label{1.1}
S_N=\sum_{n=1}^NF(\xi_n,\xi_{2n},...,\xi_{\ell n}).
\end{equation}
When $F(x_1,...,x_\ell)=\prod^\ell_{j=1}\bbI_\Gam(x_j)$, where $\bbI_\Gam(x)
=1$ if $x\in\Gam$ and $=0$ if $x\not\in\Gam$, then $S_N$ counts the number of
events $\{\xi_n,\xi_{2n},...,\xi_{\ell n}\in\Gam\},\, n\leq N$ representing
multiple returns to $\Gam$. The name "nonconventional" comes from \cite{Fu}
where $L^2$ ergodic theorems for such sums with $F(x_1,...,x_\ell)=\
\prod_{j=1}^\ell f_j(x_j)$ and $\{\xi_n\}$ generated by iterates of a measure
preserving transformation were studied. Recently, strong laws of large numbers
and the central limit theorem type results for sums (\ref{1.1}) were obtained
in \cite{Ki3} and \cite{KV}, respectively, while other related limit theorems
for such sums were studied in a series of papers of the second author and
co-authors. Nonconventional limit theorems is probably the first topic in
probability whose original motivation comes from the ergodic theory.

Modern proofs of local limit theorems adopt usually the approach based on the
the perturbation theory of Fourier operators which goes back to \cite{Na}
(see also \cite{Na2}).
This together with quasi-compactness (or spectral gap, non-arthmeticity)
assumptions yields appropriate estimates of characteristic functions of sums
under consideration implying
the result (see, for instance, \cite{HH}). The characteristic
functions of the sums (\ref{1.1}) cannot be studied via iterates of one Fourier
operator, and so we cannot benefit here from the full strength of the above
method though some elements of existing proofs of local limit theorems will
be employed here, as well. We observe that a local limit theorem can only be
formulated if the variance in the corresponding central limit theorem is
positive and since the latter question was not studied extensively before in
the nonconventional setup we deal with it here, as well. Usually, proofs of
local limit theorems rely, in particular, on a corresponding central limit
theorem result and we employ this argument also here but, in fact, we consider
a slightly different from \cite{KV} version of a nonconventional central
limit theorem which will be introduced in the next section. Our results can be
extended to sums of more general expressions $F(\xi_{q_1(n)},\xi_{q_2(n)},
...,\xi_{q_\ell(n)})$ where $q_1(n),...,q_\ell(n)$ are certain integer
valued functions similar to \cite{KV} but in order not to overload the
exposition we will restrict ourselves with the sums of the form (\ref{1.1}).
A somewhat
related but still different problem on almost sure local limit theorem (see,
for instance, \cite{GAS} and references there) will not be considered here
in the nonconventional setup.

This paper is organized as follows. In the next section we state our main
 results and comment on them. In Section \ref{sec3} we derive our results
 on positivity of the limiting variance while in Sections
 \ref{sec4}--\ref{sec6} our nonconventional local limit theorems will be proved.

\section{Preliminaries and main results}\label{sec2}\setcounter{equation}{0}

Our setup consists of a probability space $(\Om,\cF,P)$ together with
a stationary Markov chain $\xi_0,\xi_1,\xi_2,...$ evolving on a Polish
space $\cX$ equiped with a Borel $\sig$-algebra $\cB$. Let $P(x,\Gam)=
P\{\xi_1\in\Gam |\xi_0=x\}$ be the transition probability of the Markov
chain $\{\xi_n\}$, $\mu$ be its stationary probability so that
$\int d\mu(x)P(x,\Gam)=\mu(\Gam)$ for any $\Gam\in\cB$, and let $F=F(x_1,
...,x_\ell),\,\ell\geq 1$ be a Borel function on $\cX^\ell=\cX\times
\cdots\times\cX$ such that
\begin{equation}\label{2.1}
b^2=\int F^2(x_1,...,x_\ell)d\mu(x_1)\dots d\mu(x_\ell)<\infty.
\end{equation}
Our main goal in this paper is to derive a local limit theorem for sums
given by (\ref{1.1}) and in order to simplify formulas we assume also a
centering condition
\begin{equation}\label{2.2}
\bar F=\int F(x_1,...,x_\ell)d\mu(x_1)\dots d\mu(x_\ell)=0
\end{equation}
which is not a restriction since we always can replace $F$ by $F-\bar F$.

Introduce $\sig$-algebras $\cF_n=\sig\{\xi_j,\, j\leq n\}$ and $\cF^m=\sig\{
\xi_j,\, j\geq m\}$ and define the $\psi$-mixing (dependence) coefficient by
\begin{equation}\label{2.3}
\psi(m)=\sup\left\{\big\vert\frac {P(\Gam\cap\Psi)}{P(\Gam)P(\Psi)}-1\big\vert
:\,\Gam\in\cF_n,\,\Psi\in\cF^{n+m};\, P(\Gam),P(\Psi)>0\right\}.
\end{equation}
The process $\{\xi_n\}$ is called $\psi$-mixing if all $\psi(m)$ are finite
and $\psi(m)\to 0$ as $m\to\infty$. To avoid excessive technicalities we will
work here under
\begin{assumption}\label{ass2.1}
For some $\al>0$ and any $m\geq 1$,
\begin{equation}\label{2.4}
\psi(m)\leq\al^{-1}e^{-\al m}.
\end{equation}
\end{assumption}

The exponentially fast decay (\ref{2.4}) of the $\psi$-mixing coefficient
can be relaxed to some polynomial decay but it is known from \cite{BHK} that
in the Markov chains case any decay $\psi(m)\to 0$ as $m\to\infty$ yields
already an exponentially fast decay of $\psi(m)$. It is known also (see,
for instance, \cite{Br},\, Ch. 7 and 21) that (\ref{2.4}) holds true if
$\{\xi_n\}$ is a finite state irreducible and aperiodic Markov chain and,
in fact, \cite{BHK} provides necessary and sufficient conditions for
(\ref{2.4}) to take place. In particular, (\ref{2.4}) will be satisfied if
there exists a positive integer $n_0$, a probability measure $\eta$ on $\cX$
and a number $\gam\in(0,1]$ such that for each $x\in\cX$ and any Borel set
$\Gam\subset\cX$,
\begin{equation}\label{2.5}
\gam^{-1}\eta(\Gam)\geq P(n_0,x,\Gam)\geq\gam\eta(\Gam)
\end{equation}
where $P(k,x,\cdot)$ is the $k$-step transition probability of the Markov
chain $\{\xi_n\}$. Employing the technique from \cite{KV} we can
obtain our results under weaker mixing conditions on expense of other
 assumptions on the function $F$. Observe that the right hand side of
 (\ref{2.5}) implies also the geometric ergodicity condition
 \begin{equation}\label{2.5+}
 \|P(n,x,\cdot)-\mu\|\leq\be^{-1}e^{-\be n},\,\be >0,
 \end{equation}
 where $\|\cdot\|$ is the total variation norm, and under (\ref{2.5+}) the
 results of the present paper can be obtained for any initial distribution
 of the Markov chain $\{\xi_n\}$ and not only for the stationary one.

 As usual, our local limit theorem will rely on a version of a nonconventional
 central limit theorem which will be presented in a more general form than
 needed here and under slightly different assumptions than in \cite{KV}. It
 will be convenient to represent the function $F=F(x_1,...,x_\ell)$ in the
 form
 \begin{equation}\label{2.6}
 F=F_1(x_1)+\cdots +F_\ell(x_1,...,x_\ell)
 \end{equation}
 where
 \begin{equation}\label{2.7}
 F_\ell=F(x_1,...,x_\ell)-\int F(x_1,...,x_\ell)d\mu(x_\ell)
 \end{equation}
 and for $i<\ell$,
 \begin{equation}\label{2.8}
 F_i(x_1,...,x_i)=\int F(x_1,...,x_\ell)d\mu(x_{i+1})...d\mu(x_\ell)-
 \int F(x_1,...,x_\ell)d\mu(x_i)...d\mu(x_\ell)
 \end{equation}
 which ensures, in particular, that for all $x_1,...,x_{i-1}\in\cX$,
 \begin{equation}\label{2.9}
 \int F_i(x_1,...,x_{i-1},x_i)d\mu(x_i)=0.
 \end{equation}
 Now we write
 \begin{equation}\label{2.10}
 S_N(t)=\sum_{i=1}^\ell S_{i,N}(t)\,\,\,\mbox{where}\,\,\, S_{i,N}(t)=
 \sum_{n=1}^{[Nt]}F_i(\xi_n,\xi_{2n},...,\xi_{in})
 \end{equation}
 and we abbreviate $S_N=S_N(1)$ and $S_{i,N}=S_{i,N}(1)$.

 \begin{theorem}\label{thm2.2} Suppose that (\ref{2.1}), (\ref{2.2}) and
 Assumption \ref{ass2.1} hold true. Then the $\ell$-dimensional process
 $\{ N^{-1/2}S_{j,N}(t/j):\, 1\leq j\leq\ell\}$ converges in distribution
 to a Gaussian process $\{\zeta_j(t):\, 1\leq j\leq\ell\}$ with stationary
 independent increments, zero means and covariances having the form
 $E(\zeta_i(s)\zeta_j(t))=D_{i,j}\min(s,t),\, i,j=1,...,\ell$. The process
 $N^{-1/2}S_N(\cdot)$ itself converges in distribution to the Gaussian
 process $\zeta(\cdot)$ having a representation in the form
 \[
 \zeta(t)=\sum_{i=1}^\ell\zeta_i(it)
 \]
 which may have dependent increments.
 \end{theorem}

 In spite of different assumptions on the function $F$ here in comparison
 to \cite{KV} this theorem follows in the same way as the main result of
 \cite{KV} which will be explained at the beginning of Section \ref{sec3}.

 A local limit theorem can only be meaningful if the variance of a limiting
 Gaussian distribution is strictly positive. This question was addressed in
 \cite{Ki2} only in a very particular case (which corresponds to $\ell=1$
 in our setup) and it was    not dealt with at all in \cite{KV}. Here we will
 establish some sufficient conditions for positivity of the limiting variance
 in our Markov chains setup.

 \begin{theorem}\label{thm2.4-} (i) Set $\sig^2_N=\mbox{var} S_N=E(S_N-ES_N)^2$
 and suppose that (\ref{2.4}) holds true. Then the limit
 \begin{equation}\label{2.14}
 \lim_{N\to\infty}\frac 1N\sig_N^2=\sig^2=\lim_{N\to\infty}\frac 1NES_N^2
 \end{equation}
 exists.

 (ii) Let $\{\xi_n^{(1)}\}$, $\{\xi_n^{(2)}\}$, ...,$\{\xi_n^{(\ell)}\}$ be
 $\ell$ independent copies of the stationary Markov chain $\{\xi_n\}$ with
 the initial distribution $\mu$ and set $U_{\ell,N}=\sum_{n=1}^N
 F_\ell(\xi_n^{(1)},\xi_{2n}^{(2)},...,\xi_{\ell n}^{(\ell)})$. Then
 \begin{equation}\label{2.11+}
 s^2_\ell=\lim_{N\to\infty}\frac 1N\mbox{var}U_{\ell,N}=\lim_{N\to\infty}
 \frac 1N EU_{\ell,N}^2
 \end{equation}
 exists and
 \begin{equation}\label{2.11++}
 \sig^2\geq\frac 1{2\ell}s_\ell^2.
 \end{equation}
 \end{theorem}

 Of course, if $F_\ell=0$ $\mu^\ell$-a.e. then $s_\ell=0$ but in this
 case $F$ depends essentially only on $\ell-1$ variables and we can apply
 the same arguments with $\ell-1$ in place of $\ell$. Next, observe that
 specifying a bit our estimates it is possible to improve the right hand side
 of (\ref{2.11++}) to $\frac 1\ell s^2_\ell$. Nevertheless, the main purpose
 of Theorem \ref{thm2.4-}(ii) is to obtain a sufficient condition for
 positivity of $\sig^2$ relying on the well known conditions of positivity of
 the limiting variance of the stationary Markov chain $\Xi_n^{(\ell)}=
 (\xi_n^{(1)},\xi_{2n}^{(2)},...,\xi_{\ell n}^{(\ell)})$ which is $\psi$-mixing
 with the same $\psi(n)$-coefficient as $\{\xi_n\}$ itself. Namely, $s_\ell^2>0$
 unless $F_\ell$ can be represented in the form $F_\ell(\Xi_1^{(\ell)}(\om))=
 g(\te\om)-g(\om)$ where $\te$ is the $P$-preserving paths' shift
 transformation for the stationary process $\{\Xi_n^{(\ell)}\}$ so that
 $\Xi_n^{(\ell)}=\Xi_0^{(\ell)}(\te^n\om)$ (see Chapter 18 in \cite{IL}).
 We observe also that the assertion of Theorem \ref{thm2.4-} remain valid
 for any stationary (and not only Markov) process with a $\psi$-mixing
 coefficient decaying much slower than in (\ref{2.4}) (see the proof in the
  next section). Furthermore, these assertions hold also true for sufficiently
   fast $\psi$-mixing dynamical systems (such as subshifts of finite type,
   Anosov diffeomorphisms, expanding transformations etc.) where in place of
   independent copies of the process $\{\xi_n\}$ we should consider the product
   of $\ell$ copies of the corresponding dynamical system.

 Introduce correlation coefficients
 \begin{equation}\label{2.11}
 \rho_k=\sup_{n\geq 1}\{\mbox{corr}(X,Y):\, X\in L_2(\Om,\cF_n,P),\,
 Y\in L_2(\Om,\cF^{n+k},P)\}=\| Q^k\|_{L^0_2(\cX,\mu)},
 \end{equation}
 where $Qg(x)=\int_\cX P(x,dy)g(y)$ is the transition operator and
 $L^0_2(\cX,\mu)=\{ g\in L_2(\cX,\cB,\mu):\,\int gd\mu=0\}$ . We consider also
 contraction coefficients
 \begin{equation}\label{2.12}
 \del_k=\sup_{x,y\in\cX,\Gam\in\cB}|P(k,x,\Gam)-P(k,y,\Gam)|.
 \end{equation}
 In some of the assertions we will assume that
 \begin{equation}\label{2.13}
 \rho_\ell<1.
 \end{equation}
 It was indicated to us by M. Peligrad that the arguments in the proof
 of Lemma 4.1 of \cite{SV} yield that $\rho_\ell\leq\sqrt\del_\ell$, and so if
 \begin{equation}\label{2.13+}
 \del_\ell<1
 \end{equation}
 then (\ref{2.13}) holds true.

 Next, we provide additional more specific sufficient conditions for
 the positivity of $\sig^2$.
 \begin{theorem}\label{thm2.4}
(i) Let (\ref{2.4}) and (\ref{2.13}) hold true. If
 $\sig^2=0$ then $F_\ell=0$ $\mu^\ell$-a.e.

 (ii) Suppose that (\ref{2.4}) and (\ref{2.13})
 hold true and assume that $P(x,\cdot)$ is absolutely continuous with respect
 to $\mu$ for all $x\in\mbox{supp}\mu$. Then $\sig^2=0$ if and only if
 $F=0$ $\mu^\ell=\mu\times\cdots\times\mu$-almost everywhere (a.e.);

 (iii) Suppose that the conditions of (ii) hold true except that in place of
 the absolute continuity requirement there we assume here that $F$ is
 continuous on (supp$\mu)^\ell=\mbox{supp}\mu\times\cdots\times\mbox{supp}\mu$.
 Then $\sig^2=0$ if and only if $F=0$ identically on (supp$\mu)^\ell$.
 \end{theorem}

 \begin{remark}\label{rem2.4+} In principle, we could exclude $F_\ell=0$
 $\mu^\ell$-a.e.
  from the beginning since if this happens then $F$ can be viewed
 as a function of $\ell-1$ variables and we could study then the problem
 with $\ell-1$ in place of $\ell$. Thus, from the beginning we could assume
 that $\ell$ is the maximal number such that $F_\ell=0$ $\mu^\ell$-a.e. does
 not happen.
 \end{remark}

 Observe that when $\xi_0,\xi_1,\xi_2,...$ are independent identically
 distributed (i.i.d.) random variables with values in $\cX$ then $\rho_k=
 \del_k=\psi(k)=0$ for all $k\geq 1$, and so (\ref{2.4}), (\ref{2.13})
 and (\ref{2.13+}) hold true. The inequality (\ref{2.13}) is known to
 hold true, for instance, when the Markov chain $\{\xi_n\}$ is reversible and
 geometrically ergodic (see \cite{RR}). Some spectral conditions for
 (\ref{2.13}) to hold true can be found in \cite{Br}. The validity of
 (\ref{2.13+}) can be ensured assuming only the right hand side of (\ref{2.5})
  with  $n_0=\ell$. When $\cX$ is a finite state space
 conditions for (\ref{2.13}) and (\ref{2.13+}) to be valid can be
 easily written (see, for instance, \cite{Do}, \cite{Br} and referencies
 there). We will need (\ref{2.13}) (or (\ref{2.13+}))
 in order to rely on lower bounds for variances of sums of (different)
 functions of Markov chains from \cite{Pe} while the corresponding
 lower bound from \cite{SV} can also be used if (\ref{2.13+}) holds true.

 As in many expositions of the (conventional) local limit theorem we
 distinguish between a lattice and a nonlattice cases which take in our
 circumstances the following form. For any $x_1,...,x_{\ell-1}\in\cX$ set
 \begin{equation}\label{2.15}
 A_{x_1,...,x_{\ell-1}}=\{ h\geq 0:\, F(x_1,...,x_{\ell-1},x)
 \in\{ kh:\, k\in\bbZ\}\,\,\mbox{for}\,\,\mu-\mbox{almost all}\,\, x\in\cX\}
 \end{equation}
 and
  \begin{eqnarray}\label{2.15+}
 &B_{x_1,...,x_{\ell-1}}=\{ h\geq 0:\, F(x_1,...,x_{\ell-1},x)-F(x_1,...,
 x_{\ell-1},y)\\
 &\in\{ kh:\, k\in\bbZ\}\,\,\mbox{for}\,\,\mu^2-\mbox{almost all}\,\, (x,y)
 \in\cX^2\}.\nonumber
 \end{eqnarray}
 If $B_{x_1,...,x_{\ell-1}}\ne\emptyset$ we define
 \begin{equation}\label{2.16}
 h(x_1,...,x_{\ell-1})=\sup\{ h:\, h\in B_{x_1,...,x_{\ell-1}}\}.
 \end{equation}
 We call the case a lattice one if there exists $h>0$ such that
 \begin{equation}\label{2.17}
 h(x_1,...,x_{\ell-1})=h\in A_{x_1,...,x_{\ell-1}} \quad\,\,\,{for}\,\,\,
 \mu^{\ell-1}-\mbox{almost all}\,\,\, (x_1,...,x_{\ell-1}).
 \end{equation}
 If
 \begin{equation}\label{2.18}
 \mu^{\ell-1}\{(x_1,...,x_{\ell-1}):\, B_{x_1,...,x_{\ell-1}}=\emptyset\}>0
 \end{equation}
 then we call the case a non-lattice one. Observe, that there are other
 cases beyond what we designated as a lattice and a non-lattice case. For
 instance, $h(x_1,...,x_{\ell-1})$ may take on several (or countably many,
 or continuum) values on subsets of $\cX^{\ell-1}$ having positive (or zero
 in the continuum values case) $\mu^{\ell-1}$-measure. Some of these cases
 can be treated too but in order not to overload the exposition we do not
 consider them here.

 Now we can state our nonconventional local limit theorems which are considered
 both in a lattice and a non-lattice cases.

 \begin{theorem}\label{thm2.6} Suppose that $\xi_0,\xi_1,\xi_2,...$ is a
 Markov chain having transition probabilities satisfying (\ref{2.5})
  with $n_0=\ell$. Assume (\ref{2.1}) and (\ref{2.2}) with $b\ne 0$ and that
  the equality $F_\ell=0\,\,\,\mu^\ell$-a.e. does not hold true. Then
  $\sig^2$ in (\ref{2.14}) is positive and for any real continuous function
  $g$ on $\bbR$ with a compact support
 \begin{equation}\label{2.19}
 \lim_{N\to\infty}\sup_{u}|\sig\sqrt {2\pi N}Eg(S_N-u)-e^{-\frac {u^2}
 {2N\sig^2}}\int gd\cL_h|=0
 \end{equation}
 where the supremum is taken over $u\in\{ kh:\, k\in\bbZ\}$ in the lattice case
 (\ref{2.17}) and over all real $u$ in the non-lattice case. In the latter case
 we set $h=0$ and $\cL_0$ is then the Lebesgue measure on $\bbR$ while in the
 lattice case $\cL_h,\, h>0$ is the measure
 assigning mass $h$ to each point $kh,\, k\in\bbZ$.
 \end{theorem}

 The proof of Theorem \ref{thm2.6} will rely on estimates of the characteristic
 function $\vf_N(\te)$ of $S_N$ of the form $e^{-qN}$ for $\te$
 belonging to a compact set disjoint from 0 and of the form
  $e^{-rN\te^2}$ for $\te$ close to 0. These are obtained employing
  some large deviations results from \cite{DV}, \cite{Bol} and \cite{KM}.
  From these estimates Theorem \ref{thm2.6} will follow essentially in a
  standard way. When $\cX$ is a finite state space Theorem \ref{thm2.6}
  requires that the transition matrix of the Markov chain $\{\xi_n\}$ consists
  of positive entries only. Actually, in this case we can obtain the
  nonconventional local limit theorem in a bit more general situation.

  \begin{theorem}\label{thm2.7} Let $\{\xi_n\}$ be a Markov chain with a
  finite state space $\cX$ and a matrix of transition probabilities
  $\Pi=(p_{ij})_{i,j\in\cX}$ such that for some $k$ the matrix
  $\Pi^{k}=(p^{(k)}_{ij})_{i,j\in\cX}$ consists of all positive
  entries and
  \begin{equation}\label{2.20}
  \min_{i,j}\sum_k\min(p^{(\ell)}_{ik},p_{jk}^{(\ell)})>0\,\,\mbox{and}\,\,
  \min_{i,j}\sum_k\min(p^{(\ell)}_{ki},p_{kj}^{(\ell)})>0.
  \end{equation}
  Assume (\ref{2.2}) and that $F_\ell$ does not equal zero
  identically. Then both in the lattice and the non-lattice cases
  the local limit theorem in the form (\ref{2.19}) takes place.
  \end{theorem}
  We observe that (\ref{2.13+}) implies that $\del_\ell<1$ and the latter
  is equivalent to the first inequality in (\ref{2.20}) (see \cite{Do}).
  The first and the second inequalities in (\ref{2.20}) mean that the matrices
  $\Pi^\ell(\Pi^{\ell})^*$ and $(\Pi^{\ell})^*\Pi^\ell$, respectively, have all
   positive entries. The conditions of Theorem \ref{thm2.7} concerning $\Pi$
   hold true if, for instance, $\Pi$ is self adjoint and $\Pi^2$ has all
   positive entries.

\section{Positivity of variance}\label{sec3}\setcounter{equation}{0}

In this section we explain first why the proof of Theorem \ref{thm2.2}
goes through in the same way as in \cite{KV} and after that we derive
Theorem \ref{thm2.4}.

We will rely on the following result whose proof goes through in the same
way as in Lemma 3.1 from \cite{Ki1}. Let $h=h(y,z)$, $y\in\cX^k,\,
z\in\cX^l$ be a Borel measurable function on $\cX^{k+l}$. Let
$X=(X_1,...,X_k):\,\Om\to\cX^k$ and $Y=(Y_1,...,Y_l):\,\Om\to\cX^l$ be
sequences of $\cX$-valued random variables (random "vectors") such that
$X$ is $\cF_n$-measurable and $Y$ is $\cF^{n+m}$-measurable. Suppose
that $E|h(X,Y)|<\infty$ and $r(x)=E|h(x,Y)|<\infty$ for any $x\in\cX^k$.
Then with probability one,
\begin{equation}\label{3.1}
|E(h(X,Y)|\cF_n)-g(X)|\leq\psi(m)r(X)
\end{equation}
where $g(x)=Eh(x,Y)$.

Now observe that special H\" older continuity assumptions imposed on $F$ in
 \cite{KV} were only needed because we did not assume there measurability
 of $\xi_n$ with respect to $\cF_n$ but required only some estimates on the
 error when replacing $\xi_n$ by $E(\xi_n|\cF_{n+m})$. To use this we had to
 estimate what happens when $\xi_n$ is replaced by $E(\xi_n|\cF_{n+m})$ inside
 $F$ which could be done only relying on H\" older continuity (or better
 regularity) of $F$. In the present paper $\cF_n$ is generated by $\xi_0,\xi_1,
 ...,\xi_n$, i.e. $\xi_n$ is $\cF_n$-measurable, and we can apply (\ref{3.1})
 directly to $F$ without any need of H\" older continuity as exhibited below.

 We will show next how to obtain variance and covariance estimates needed to
 employ the method from \cite{KV} to the setup of Theorem \ref{thm2.2} using
 directly (\ref{3.1}) without any need of regularity assumptions on $F$.
 Define inductively for $i=\ell-1,\ell-2,...,0$,
 \begin{equation}\label{3.2}
 G_i^2(x_1,...,x_i)=EG^2_{i+1}(x_1,...,x_i,\xi_1)=\int G^2_{i+1}(x_1,...,x_i,y)
 d\mu(y)
 \end{equation}
 where $G_\ell^2=F^2$ and $G^2_0=EG^2_1(\xi_1)$. Then by (\ref{3.1}) considered
 with $X=(\xi_n,\xi_{2n},...,\xi_{(i-1)n})$ and $Y=\xi_{in}$, $i=\ell,\ell-1,
 ...,1$ we obtain
 \begin{eqnarray}\label{3.3}
 &EF^2(\xi_n,\xi_{2n},...,\xi_{\ell n})\leq (1+\psi(n))EG_{\ell-1}^2(\xi_n,...,
 \xi_{(\ell-1)n})\\
 &\leq\cdots\leq (1+\psi(n))^\ell G^2_0=(1+\psi(n))^\ell b^2\nonumber
 \end{eqnarray}
 where $b^2$ is given by (\ref{2.1}). Next, observe that by the Cauchy-Schwarz
 inequality
 \begin{eqnarray}\label{3.4}
 &F^2_i(x_1,...,x_i)\leq 2(\int F(x_1,...,x_\ell)d\mu(x_{i+1})...
 d\mu(x_\ell))^2\\
 &+2(\int F(x_1,...,x_\ell)d\mu(x_{i})...d\mu(x_\ell))^2
 \leq 2\int F^2(x_1,...,x_\ell)d\mu(x_{i+1})...d\mu(x_\ell)\nonumber\\
 &+2\int F^2(x_1,...,x_\ell)d\mu(x_{i})...d\mu(x_\ell)
 =2G^2_i(x_1,...,x_i)+2G^2_{i-1}(x_1,...,x_{i-1}),\nonumber
 \end{eqnarray}
 and so, as in (\ref{3.3}),
 \begin{equation}\label{3.5}
 EF^2_i(\xi_n,...,\xi_{in})\leq 2(1+\psi(n))^{i-1}(2+\psi(n))b^2.
 \end{equation}

 Now we estimate covariances for $j\geq i$ and $m>k$ which by (\ref{3.1})
 considered with $X=\xi_k,...,\xi_{ik};\xi_m,...,\xi_{(j-1)m}$ and $Y=\xi_{jm}$
 together with (\ref{2.9}) yields
 \begin{eqnarray}\label{3.6}
 &A_{i,j,k,m}=|E(F_i(\xi_k,...,\xi_{ik})F_j(\xi_m,...,\xi_{jm}))|\\
 &\leq\psi(\min(m,jm-ik))ER(\xi_k,...,\xi_{ik};\xi_m,...,\xi_{(j-1)m})
 \nonumber\end{eqnarray}
 where
 \begin{eqnarray*}
 &R(x_1,...,x_i;y_1,...,y_{j-1})=|F_i(x_1,...,x_i)|E|F_j(y_1,...,y_{j-1},\xi_1)|\\
 &\leq\frac 12F^2_i(x_1,...,x_i)+\frac 12EF^2_j(y_1,...,y_{j-1},\xi_1).
 \end{eqnarray*}
 Hence, by (\ref{3.5}),
 \begin{eqnarray}\label{3.7}
 &ER(\xi_k,...,\xi_{ik};\xi_m,...,\xi_{(j-1)m})\leq \big((1+\psi(k))^{i-1}
 (2+\psi(k))\\
 &+(1+\psi(m))^{j-1}(2+\psi(m))\big)b^2.\nonumber
 \end{eqnarray}
 Observe also that by (\ref{2.9}), (\ref{3.1}) and the Cauchy-Schwarz
 inequality for all $i=1,...,\ell$,
 \begin{equation}\label{3.8}
 |EF_i(\xi_n,...,\xi_{in})|\leq\psi(n)EQ(\xi_n,...,\xi_{(i-1)n})
 \end{equation}
 where
 \[
 Q(x_1,...,x_{i-1})=E|F_i(x_1,...,x_{i-1},\xi_1)|\leq (EF_i^2(x_1,...,x_{i-1},
 \xi_1))^{1/2}.
 \]
 Again, by (\ref{3.5}),
 \begin{equation}\label{3.9}
 EQ(\xi_n,...,\xi_{(i-1)n})\leq\sqrt 2(1+\psi(n))^{(i-1)/2}(2+\psi(n))^{1/2}|b|.
 \end{equation}
 The estimates (\ref{3.4})--(\ref{3.9}) suffice in order to employ martingale
 approximations and to obtain limiting covariances as in \cite{KV} which
 yields Theorem \ref{thm2.2} in the same way as there. \qed

 Next, we deal with Theorem \ref{thm2.4-}. First observe that the proof of
 existence of the limit (\ref{2.14}) is the same as in \cite{KV} relying
 on the variance and covariance estimates above, and so we concentrate on
 the proof of the assertion (ii) of Theorem \ref{thm2.4-}. Set $N_\ell=
 [N(1-\frac 1{2\ell})]+1$,
 \begin{equation}\label{3.10}
 S^{(1)}_{\ell,N}=\sum_{n=1}^{N_\ell-1}F_\ell(\xi_n,...,\xi_{\ell n})\,\,\,
 \mbox{and}\,\,\, S^{(2)}_{\ell,N}=S_{\ell,N}-S_{\ell,N}^{(1)}.
 \end{equation}
 Thus,
 \begin{eqnarray}\label{3.11}
 &\sig^2_N=\mbox{var}S_N=\mbox{var}(\sum_{i=1}^{\ell-1}S_{i,N}+S^{(1)}_{\ell,N})
 +\mbox{var}S^{(2)}_{\ell,N}\\
 &+2\sum_{i=1}^{\ell-1}\sum_{k=1}^N\sum_{m=N_\ell}^N B_{i,k,m}+
 2\sum_{k=1}^{N_\ell-1}\sum_{m=N_\ell}^NB_{\ell,k.m}\nonumber
 \end{eqnarray}
 where for $i=1,2,...,\ell,$
 \[
 B_{i,k,m}=E\big((F_i(\xi_k,...,\xi_{ik})-EF_i(\xi_k,...,\xi_{ik}))
 (F_\ell(\xi_m,...,\xi_{\ell m})-EF_\ell(\xi_m,...,\xi_{\ell m}))\big).
 \]

 In order to prove Theorem \ref{thm2.4-}(ii) we will show first
 that the multiple sums in (\ref{3.11}) are bounded after which it will
 remain only to obtain that
 \begin{equation}\label{3.12}
 \liminf_{N\to\infty}\frac 1N\mbox{var} S^{(2)}_{\ell,N}>0.
 \end{equation}
 Clearly,
 \begin{equation}\label{3.13}
 |B_{i,k,m}|\leq A_{i,\ell,k,m}+|EF_i(\xi_k,...,\xi_{ik})||EF_\ell(\xi_m,
 ...,\xi_{\ell m})|.
 \end{equation}
 Observe that
 \begin{equation}\label{3.14}
 \ell N_\ell-(\ell-1)N\geq N/2,
 \end{equation}
 and so when $m\geq N_\ell,\, i\leq\ell-1$ and $k\leq N$ we can estimate the
 right hand side of (\ref{3.13}) by means of (\ref{2.4}) and
 (\ref{3.6})--(\ref{3.9}) which yields that two double sums in (\ref{3.11}) are
 bounded by constants independent of $N$.

 Next, we study
 \begin{eqnarray}\label{1+}
 &B_{\ell,m,n}=E\big(F_\ell(\xi_m,\xi_{2m},...,\xi_{\ell m})
 F_\ell(\xi_n,\xi_{2n},...,\xi_{\ell n})\big)\\
 &-EF_\ell(\xi_m,\xi_{2m},...,\xi_{\ell m})
 EF_\ell(\xi_n,\xi_{2n},...,\xi_{\ell n})\nonumber
 \end{eqnarray}
 for $m,n\geq N_\ell$ in more details. The product of expectations in
 (\ref{1+}) is small by (\ref{3.8}) and (\ref{3.9}), and so we have to
 deal mainly with the expectation of the product here. Set
 \[
 H^{(\ell)}_{m,n}(x_1,...,x_\ell;\, y_1,...,y_\ell)=F_\ell(x_1,...,x_\ell)
 F_\ell(y_1,...,y_\ell)
 \]
 and define recursively for $j=\ell,\ell-1,...,1$,
 \[
 H^{(j-1)}_{m,n}(x_1,...,x_{j-1};\, y_1,...,y_{j-1})=\int H^{(j)}_{m,n}(x_1,...,
 x_j;\, y_1,...,y_j)d\mu_{j|n-m|}(x_j,y_j)
 \]
 where $\mu_k$ is the distribution of the pair $(\xi_0,\xi_k)$.

 Next, observe that
 \begin{equation}\label{2+}
 H^{(0)}_{m,n}=E\big(F_\ell(\xi^{(1)}_m,\xi_{2m}^{(2)},...,
 \xi^{(\ell)}_{\ell m})F_\ell(\xi^{(1)}_n,\xi_{2n}^{(2)},...,
 \xi^{(\ell)}_{\ell n})\big)
 \end{equation}
 and that $EF_\ell(\xi^{(1)}_m,\xi_{2m}^{(2)},...,
 \xi^{(\ell)}_{\ell m})=0$ by (\ref{2.2}). Applying (\ref{3.1}) with
 $X=(\xi_m,...,\xi_{(j-1)m};\,\xi_n,...,\xi_{(j-1)n})$ and
 $Y=(\xi_{jm},\xi_{jn})$ we obtain that for any $m,n\geq N_\ell$,
 \begin{eqnarray}\label{3+}
 &\big\vert EH^{(j)}_{m,n}(\xi_m,\xi_{2m},...,\xi_{jm};\,\xi_n,\xi_{2n},
 ...,\xi_{jn})\\
 &-EH^{(j-1)}_{m,n}(\xi_m,\xi_{2m},...,\xi_{(j-1)m};\,\xi_n,\xi_{2n},
 ...,\xi_{(j-1)n})\big\vert\nonumber\\
 &\leq\psi([\frac N2])E\hat H^{(j-1)}_{m,n}(\xi_m,\xi_{2m},...,\xi_{(j-1)m};
 \,\xi_n,\xi_{2n},...,\xi_{(j-1)n})\nonumber
 \end{eqnarray}
 where
 \[
 \hat H^{(j-1)}_{m,n}(x_1,...,x_{j-1}; \, y_1,...,y_{j-1})=E|H^{(j)}_{m,n}(x_1,
 ...,x_{j-1},\xi_{jm};\, y_1,...,y_{j-1},\xi_{jn})|.
 \]
 Since
 \begin{eqnarray*}
 &H^{(\ell)}_{m,n}(x_1,...,x_{\ell}; \, y_1,...,y_{\ell})=
 |F_\ell(x_1,...,x_\ell)||F_\ell(y_1,...,y_\ell)|\\
 &\leq\frac 12F_\ell^2(x_1,...,x_\ell)+\frac 12F_\ell^2(y_1,...,y_\ell)
 \end{eqnarray*}
 we obtain recursively from the above that for $j=\ell-1,\ell-2,...,0$,
 \begin{equation}\label{4+}
 \hat H^{(j)}_{m,n}(x_1,...,x_{j}; \, y_1,...,y_{j})\leq 2G^2_j(x_1,
 ...,x_j)+2G^2_j(y_1,...,y_j)
 \end{equation}
 with $G^2_j$ defined by (\ref{3.2}). In the same way as in (\ref{3.3}) we see
 that
 \begin{equation}\label{5+}
 EG_j^2(\xi_n,...,\xi_{jn})\leq (1+\psi(n))^jb^2.
 \end{equation}
 Now, combining (\ref{3+})--(\ref{5+}) we obtain that
\begin{equation}\label{6+}
\big\vert E\big(F_\ell(\xi_m,\xi_{2m},...,\xi_{\ell m})F_\ell(\xi_n,\xi_{2n},
...,\xi_{\ell n})\big)-H^{(0)}_{m,n}\big\vert\leq\psi([\frac N2])\ell(1+
\psi(N_\ell))^\ell b^2.
\end{equation}

 Set
 \[
 U^{(2)}_{\ell,N}=\sum^N_{n=N_\ell}F_\ell(\xi_n^{(1)},\xi_{2n}^{(2)},...,
 \xi_{\ell n}^{(\ell)}).
 \]
 Then by (\ref{3.8}), (\ref{3.9}), (\ref{1+}), (\ref{2+}) and (\ref{6+}),
 \begin{equation}\label{7+}
 |\mbox{var}S^{(2)}_{\ell,N}-\mbox{var}U^{(2)}_{\ell,N}|\leq\frac {N^2}{4\ell^2}
 \psi([\frac N2])(2+\psi([\frac N2]))^\ell b^2(\ell+2\psi([\frac N2])).
 \end{equation}
 Next, observe that since $\Xi^{(\ell)}_n=(\xi_n^{(1)},\xi_{2n}^{(2)},...,
 \xi_{\ell n}^{(\ell)})$ is a stationary process then
 \begin{equation}\label{8+}
 \mbox{var}U^{(2)}_{\ell,N}=\mbox{var}U_{\ell,N-N_\ell}.
 \end{equation}
 It is a standard fact that under (\ref{2.4}) (and even much slower
 $\psi$-mixing) the limit (\ref{2.11+}) exists (see, for instance, \cite{IL}
 or \cite{Br}), and so
 \begin{equation}\label{9+}
 \lim_{N\to\infty}\frac 1N\mbox{var}U^{(2)}_{\ell,N}=\frac 1{2\ell}s^2_\ell.
 \end{equation}
 Now, (\ref{3.12}) follows from (\ref{7+}) and (\ref{9+}) so that
 (\ref{2.11++}) holds true in view of (\ref{3.11}) and (\ref{3.13}),
 completing the proof of Theorem \ref{thm2.4-}.     \qed

 Next we prove Theorem \ref{thm2.4}.
 In order to obtain (\ref{3.12}) we set for each $\bar x=(x_{N_\ell},
 x_{N_\ell +1},...,x_{(\ell-1)N})$,
 \[
 S^{(2)}_{\ell,N}(\bar x)=\sum^N_{n=N_\ell}F_\ell(x_n,x_{2n},...,x_{(\ell-1)n},
 \xi_{\ell n})
 \]
 where in the case $\ell=1$ we take $S^{(2)}_{\ell,N}(\bar x)=
 S^{(2)}_{\ell,N}$. Put also
 $\sig^2_{N_\ell,N}(\bar x)=$var$S^{(2)}_{\ell,N}(\bar x)$. Then applying
 (\ref{3.1}) with $X=\bar\xi=(\xi_{N_\ell},\xi_{N_\ell+1},...,\xi_{(\ell-1)N})$
 and $Y=(\xi_{\ell N_\ell},...,\xi_{\ell N})$ we obtain taking into account
 (\ref{3.5}) and (\ref{3.14}) that
 \begin{equation}\label{3.15}
 |E\sig^2_{N_\ell,N}(\bar\xi)-\mbox{var}S^{(2)}_{\ell,N}|\leq
 2N^2(2+\psi(1))^{\ell+1}\psi([N/2])b^2.
 \end{equation}

 Next, we apply to $S^{(2)}_{\ell,N}(\bar x)$ Proposition 13 from
 \cite{Pe} (if we assume (\ref{2.13+}) then we can rely on Proposition 4.1
  from \cite{SV}) in order to conclude that
 \begin{equation}\label{3.16}
 \mbox{var}S^{(2)}_{\ell,N}(\bar x)\geq c_\ell\sum_{n=N_\ell}^N\mbox{var}
 F_\ell(x_n,...,x_{(\ell-1)n},\xi_{\ell n})
 \end{equation}
 where $c_\ell=\frac {1-\rho_\ell}{1+\rho_\ell}$ under (\ref{2.13})
 (and $c_\ell=\frac {1-\del_\ell}4$ under (\ref{2.13+})).
 Observe that by (\ref{2.9}) and the stationarity, for any $y_1,...,y_{\ell-1}$,
 \begin{eqnarray*}
 &q(y_1,...,y_{\ell-1})=\mbox{var}F_\ell(y_1,...,y_{\ell-1},\xi_{\ell n})\\
 &=\mbox{var}F_\ell(y_1,...,y_{\ell-1},\xi_1)=EF^2_\ell(y_1,...,y_{\ell-1},
 \xi_1).
 \end{eqnarray*}
 Applying $\ell-1$ times (\ref{3.1}) and using (\ref{3.6}) we obtain easily
  that
 \begin{eqnarray}\label{3.17}
 &|Eq(\xi_n,\xi_{2n},...,\xi_{(\ell-1)n})-\int q(x_1,...,x_{\ell-1})d\mu(x_1)
 ... d\mu(x_{\ell-1})|\\
 &\leq 2(\ell-1)\psi(n)(2+\psi(n))^{\ell+1}b^2.\nonumber
 \end{eqnarray}
 This together with (\ref{2.4}), (\ref{3.11}), (\ref{3.13}), (\ref{3.15})
 and (\ref{3.16}) yields that
 \begin{equation}\label{3.18}
 \liminf_{N\to\infty}\frac 1N\mbox{var} S^{(2)}_{\ell,N}\geq\frac {c_\ell}\ell
 \int q(x_1,...,x_{\ell-1})d\mu(x_1)... d\mu(x_{\ell-1}).
 \end{equation}
 If $\sig^2=0$ then
 \begin{equation}\label{3.19}
 \int q(x_1,...,x_{\ell-1})d\mu(x_1)... d\mu(x_{\ell-1})=\int F^2_\ell(x_1,...,
 x_\ell)d\mu(x_1)...d\mu(x_\ell)=0,
 \end{equation}
 i.e.
 \begin{equation}\label{3.20}
 F_\ell(x_1,...,x_\ell)=0\quad\mu^\ell-\mbox{a.e.}
 \end{equation}
 which proves the assertion (i) of Theorem \ref{thm2.4}.

 Next, introduce probability measures $\mu_{i,n}$ on $\cX^i$ such that for
  any Borel $\Gam\subset\cX^i$,
  \[
  \mu_{i,n}(\Gam)=P\{(\xi_n,\xi_{2n},...,\xi_{in})\in\Gam\}.
  \]
  If $\Gam=\Gam_1\times\Gam_2\times\cdots\times\Gam_i$ is a product set
  with Borel $\Gam_1,...,\Gam_i\subset\cX$ then
  \[
  \mu_{i,n}(\Gam)=\int_{\Gam_1}d\mu(x_1)\int_{\Gam_2}P(n,x_1,dx_2)...
  \int_{\Gam_i}P(n,x_{i-1},dx_i).
  \]
  Now assume that we are in the circumstances of Theorem \ref{thm2.4}(ii),
  i.e. that $P(x,\cdot)$ is absolutely continuous with respect to $\mu$ for
  all $x\in$supp$\mu$. Then by the Chapman--Kolmogorov formula $P(n,x,\cdot)$
  is also absolutely continuous with respect to $\mu$ for all $x\in$supp$\mu$.

  We claim now that $\mu_{i,n}$ is absolutely continuous with respect to
  $\mu^i=\mu\times\cdots\times\mu$ which will be proved by induction in $i$.
  For $i=1$ this is clear since $\mu_{1,n}=\mu$. Suppose this holds true
  for $i=j$ and prove it for $i=j+1$. Let $\Gam\subset\cX^{j+1}$ and
  $\mu^{j+1}(\Gam)=0$. By the Fubini theorem for $\mu^j$-almost all
  $(x_1,...,x_j)$ the set $\Gam_{x_1,...,x_j}=\{ x\in\cX:\, (x_1,...,x_j,
  x)\in\Gam\}$ satisfies $\mu(\Gam_{x_1,...,x_j})=0$. By the induction
  hypothesis $\mu_{j,n}\prec \mu^j$, and so $\mu(\Gam_{x_1,...,x_j})=0$ for
  $\mu_{j,n}$-almost all $(x_1,...,x_j)$. From the definition of measures
  $\mu_{i,n}$, the disintegration of measures theorem (see \S 3, Ch. 6 in
  \cite{Bo}) and the Markov property we conclude that
  \[
  \mu_{j+1,n}(\Gam)=\int_{\cX^j}d\mu_{j,n}(x_1,...,x_j)P(n,x_j,
  \Gam_{x_1,...,x_j})=0
  \]
  since $P(n,x,\cdot)\prec\mu$, completing the induction step.

  Now, since $\mu_{\ell,n}\prec\mu$ then $F_\ell(\xi_n,\xi_{2n},...,
  \xi_{\ell n})=0$ with probability one for all $n$, and so $S_{\ell,N}=0$
  $P$-a.s. Hence,
  \begin{equation}\label{3.21}
  S_N=\sum_{i=1}^{\ell-1}S_{i,N}\quad P-\mbox{a.s.}
  \end{equation}
  Repeating the above argument for $\ell-1$ in place of $\ell$ we obtain that
  $F_{\ell-1}(x_1,...,x_{\ell-1})=0$ $\mu^{\ell-1}$-a.e. and
  $F_{\ell-1}(\xi_n,\xi_{2n},...,\xi_{(\ell-1) n})=0$ $P$-a.s. for all $n$,
  and so $S_{\ell-1,N}=0$ $P$-a.s., as well. Continuing in the same way
  we obtain finally from (\ref{2.6}) that $F(x_1,...,x_\ell)=0$ $\mu^\ell$-a.e.
  proving the assertion (ii).

  In the assertion (iii) we assume that $F$ is continuous on $($supp$\mu)^\ell$
  and then the absolute continuity arguments above are not needed since from
  (\ref{3.20}) we obtain directly in this case that $F_\ell(x_1,...,x_\ell)=0$
  for all $x_1,...,x_\ell\in$supp$\mu$. Since
  \[
  P\{(\xi_n,...,\xi_{\ell n})\not\in(\mbox{supp}\mu)^\ell\}\leq\sum_{i=1}^\ell
  P\{\xi_{in}\not\in\mbox{supp}\mu\}=\ell\mu(\cX\setminus\mbox{supp}\mu)=0
  \]
  we obtain again that $F_\ell(\xi_n,...,\xi_{\ell n})=0$ $P$-a.s. and
  proceeding as above we conclude that $F(x_1,...,x_\ell)=0$ not only
  $\mu^\ell$-a.e. but for all $x_1,...,x_\ell\in$supp$\mu$ in view of
  continuity of $F$ here, completing the proof. \qed

 \section{Local limit theorem}\label{sec4}\setcounter{equation}{0}

In this and the next two sections we will prove Theorems \ref{thm2.6} and
\ref{thm2.7}. The arguments of the present section are rather standard
and they are valid in a more general framework while the ones in the
next sections will be heavily affected by our setups.
The proof of (\ref{2.19}) will consist mainly
of the study of the characteristic function of the sum $S_N$. First, observe
that in view of Theorem 10.7 from \cite{Bre} (see also Section 10.4 there and
 Lemma IV.5 together with arguments of Section VI.4 in \cite{HH}) it suffices
 to prove (\ref{2.19}) for all continuous complex-valued functions $g$ on
 $\bbR=(-\infty,\infty)$ such that
 \begin{equation}\label{4.1}
 \int_{-\infty}^\infty |g(x)|dx<\infty
 \end{equation}
 and its Fourier transform
 \begin{equation}\label{4.2}
 \hat g(\la)=\int_{-\infty}^\infty e^{-i\la x}g(x)dx,\,\,\la\in\bbR
 \end{equation}
 vanishes outside of a finite interval $[-L,L]$. In particular, the
 inversion formula
 \begin{equation}\label{4.3}
 g(x)=\frac 1{2\pi}\int_{-\infty}^\infty e^{i\la x}\hat g(\la)d\la
 \end{equation}
 holds true.

 Hence,
 \begin{equation}\label{4.4}
 Eg(S_N-u)=\frac 1{2\pi}\int_{-\infty}^\infty\vf_N(\la)e^{-i\la u}
 \hat g(\la)d\la
 \end{equation}
 where $\vf_N(\la)=Ee^{i\la S_N}$ is the characteristic function of $S_N$.
 Changing variables $s=\la\sig\sqrt {N}$ we obtain
 \begin{eqnarray}\label{4.5}
 &\sig\sqrt {2\pi N}Eg(S_N-u)\\
 &=\frac 1{\sqrt{2\pi}}\int_{-\infty}^\infty\vf_N(\frac s{\sig\sqrt N})
 e^{-i\frac {su}{\sig\sqrt N}}\hat g(\frac s{\sig\sqrt N})ds.\nonumber
 \end{eqnarray}
 On the other hand, from the formula for the characteristic function of
 the Gaussian distribution and the Fourier inversion formula it follows that
 \begin{equation}\label{4.6}
 e^{-\frac {u^2}{2N\sig^2}}\int gd\cL_h=\frac {\int gd\cL_h}{\sqrt {2\pi}}
 \int_{-\infty}^\infty\exp(-\frac {i\la u}{\sig\sqrt N}-\frac {\la^2}2)d\la.
 \end{equation}

 Now, in the non-lattice case we write by (\ref{4.5}) and (\ref{4.6}),
 \begin{equation}\label{4.7}
 |\sig\sqrt {2\pi N}Eg(S_N-u)-e^{-\frac {u^2}{2N\sig^2}}\int gd\cL_0|
 \leq I_1(N,T)+I_2(N,T)+I_3(N,\del)+I_4(N,\del,T)
 \end{equation}
 where
 \[
 I_1(N,T)=\frac 1{\sqrt{2\pi}}\int_{-T}^T|\vf_N(\frac \la{\sig\sqrt N})
 \hat g(\frac \la{\sig\sqrt N})-e^{-\frac {\la^2}2}\int gd\cL_0|d\la,
 \]
 \[
 I_2(T)=\frac {|\int gd\cL_0|}{\sqrt{2\pi}}\int_{|\la|>T}e^{-\frac {\la^2}2}d\la,
 \]
 \[
 I_3(N,\del)=\frac {\|\hat g\|}{\sqrt {2\pi}}\int_{L\sig\sqrt N\geq |\la|
 \geq\del\sig\sqrt N}|\vf_N(\frac \la{\sig\sqrt N})|d\la,
 \]
 \[
 I_4(N,\del,T)=\frac {\|\hat g\|}{\sqrt {2\pi}}\int_{\del\sig\sqrt N>|\la|>T}
 |\vf_N(\frac \la{\sig\sqrt N})|d\la,
 \]
 $\|\hat g\|=\sup_\la |\hat g|$ and in writing $I_3(N,\del)$ we used the
 fact that $\hat g(s)=0$ for $s\not\in[-L,L]$.

 In the lattice case (\ref{2.17}) we proceed in a slightly different way.
 First, observe that then
 \begin{equation}\label{4.8}
 \vf_N(\frac {\la}{\sig\sqrt N}+\frac {2\pi k}h)=\vf_N(\frac \la{\sig\sqrt N})
 \quad\mbox{for all}\quad k\in\bbZ.
 \end{equation}
 Set
 \[
 r(v)=\sum_{k=-\infty}^\infty\hat g(v+\frac {2\pi k}h).
 \]
 Taking into account that here $u\in\{ kh:\, k\in\bbZ\}$
  we can rewrite (\ref{4.5}) in the following way
 \begin{equation}\label{4.9}
 \sig\sqrt {2\pi N}Eg(S_N-u)=\frac 1{\sqrt {2\pi}}\int_{-\frac {\pi\sig\sqrt N}h}
 ^{\frac {\pi\sig\sqrt N}h}\vf_N(\frac \la{\sig\sqrt N})
 e^{-i\frac {\la u}{\sig\sqrt N}}r(\frac \la{\sig\sqrt N})d\la.
 \end{equation}
 This together with (\ref{4.6}) yields
 \begin{equation}\label{4.10}
 |\sig\sqrt {2\pi N}Eg(S_N-u)-e^{-\frac {u^2}{2N\sig^2}}\int gd\cL_h|
 \leq J_1(N,T)+J_2(N,T)+J_3(N,\del)+J_4(N,\del,T)
 \end{equation}
 where
 \[
 J_1(N,T)=\frac 1{\sqrt{2\pi}}\int_{-T}^T|\vf_N(\frac \la{\sig\sqrt N})
 r(\frac \la{\sig\sqrt N})-e^{-\frac {\la^2}2}\int gd\cL_h|d\la,
 \]
 \[
 J_2(T)=\frac {|\int gd\cL_h|}{\sqrt{2\pi}}\int_{|\la|>T}e^{-\frac {\la^2}2}d\la,
 \]
 \[
 J_3(N,\del)=\frac 1{\sqrt {2\pi}}\int_{\frac {\pi\sig\sqrt N}h\geq |\la|
 \geq\del\sig\sqrt N}|\vf_N(\frac \la{\sig\sqrt N})||r(\frac \la{\sig\sqrt N})|
 d\la,
 \]
 \[
 J_4(N,\del,T)=\frac 1{\sqrt {2\pi}}\int_{\del\sig\sqrt N>|\la|>T}
 |\vf_N(\frac \la{\sig\sqrt N})||r(\frac \la{\sig\sqrt N})|d\la.
 \]

 By (the central limit) Theorem \ref{thm2.2} for any $\la$,
 \begin{equation}\label{4.11}
 \lim_{N\to\infty}\vf_N(\frac \la{\sig\sqrt N})=e^{-\frac {\la^2}2}.
 \end{equation}
 Furthermore, since $g$ is continuous and integrable on $\bbR$,
 \begin{equation}\label{4.12}
 \lim_{N\to\infty}\hat g(\frac \la{\sig\sqrt N})=\int gd\cL_0.
 \end{equation}
 In the lattice case (\ref{2.17}) it follows from Ch. 10 in \cite{Bre} that
 \begin{equation}\label{4.13}
 \lim_{N\to\infty}r(\frac \la{\sig\sqrt N})=\int gd\cL_h.
 \end{equation}
 Now by (\ref{4.11})--(\ref{4.13}) and the dominated convergence theorem
 we obtain that for any $T>0$,
 \begin{equation}\label{4.14}
 \lim_{N\to\infty}I_1(N,T)=\lim_{N\to\infty}J_1(N,T)=0.
 \end{equation}
 Next, clearly,
 \begin{equation}\label{4.15}
 \lim_{T\to\infty}I_2(T)=\lim_{T\to\infty}J_2(T)=0.
 \end{equation}
 Thus, it remains to estimate $I_3(N,\del),\, I_4(N,\del,T),\, J_3(N,\del)$
 and $J_4(N,\del,T)$. In order to do this we will need two following results
 which will be proved in the next sections for the setups of Theorems \ref{thm2.6}
 and \ref{thm2.7}.

 \begin{lemma}\label{lem4.1} There exists an integer $N_0$ such that
 for all $N\geq N_0$,
 \begin{equation}\label{4.16}
 |\vf_N(\te)|\leq e^{-qN}
 \end{equation}
 which holds true in the non-lattice case for any $\te\in[\del,\del^{-1}]$
 and $\del>0$, while in the lattice case (\ref{2.17}) the inequality
 (\ref{4.16}) holds true for any $\te\in[-\frac \pi h,\frac \pi h]\setminus
 [-\del,\del]$ and $\del>0$ with some $q=q_\del>0$ depending in both cases
 only on $\del$.
 \end{lemma}
 \begin{lemma}\label{lem4.2} There exists an integer $N_0$ such that for all
 $N\geq N_0$,
 \begin{equation}\label{4.17}
 |\vf_N(\te)|\leq e^{-rN\te^2}
 \end{equation}
 which holds true whenever $|\te|\leq\del$ with some $r=r_\del>0$ depending
 only on $\del$ which is supposed to be small enough.
 \end{lemma}

 Before proving Lemmas \ref{lem4.1} and \ref{lem4.2} we will rely on them in
 order to estimate $I_3(N,\del),\, I_4(N,\del,T),\, J_3(N,\del),\,
 J_4(N,\del,T)$ and to derive (\ref{2.19}). Indeed, with $\del$ small
 enough we estimate $I_3(N,\del)$ and $J_3(N,\del)$ by Lemma \ref{lem4.1} to
 obtain
 \begin{equation}\label{4.18}
 I_3(N,\del)\leq\frac {\|\hat g\|}{\sqrt {2\pi}}(L-\del)\sig\sqrt N
 e^{-q_\del N}\to 0\quad\mbox{as}\quad N\to\infty
 \end{equation}
 and
 \begin{equation}\label{4.19}
 J_3(N,\del)\leq\frac {e^{-q_\del N}}{\sqrt {2\pi}}\sig\sqrt N
 \int^{\pi/h}_{-\pi/h}|r(s)|ds\to 0\quad\mbox{as}\quad N\to\infty.
 \end{equation}

 Next, we estimate $I_4(N,\del,T)$ and $J_4(N,\del,T)$ by Lemma \ref{lem4.2}
 to obtain
 \begin{equation}\label{4.20}
  I_4(N,\del,T)\leq \frac {\|\hat g\|}{\sqrt {2\pi}}\int_{|\la|>T}
  e^{-r\la^2/\sig^2}d\la\to 0\quad\mbox{as}\quad T\to\infty
  \end{equation}
  and
  \begin{equation}\label{4.21}
  J_4(N,\del,T)\leq \frac {R}{\sqrt {2\pi}}\int_{|\la|>T}e^{-r\la^2/\sig^2}
  d\la\to 0\quad\mbox{as}\quad T\to\infty
  \end{equation}
 where
 \[
 R=\sup_{|u|\leq\del}|r(u)|\leq hL\sup_{|v|\leq L+\del}|\hat g(v)|<\infty.
 \]
 Thus, letting first $N\to\infty$ and then $T\to\infty$ we derive (\ref{2.19})
 from (\ref{4.14}), (\ref{4.15}) and (\ref{4.18})--(\ref{4.21}).  \qed

\section{Positive transition densities case}\label{sec5}
\setcounter{equation}{0}

 First, as it is easy to understand, (\ref{2.5}) implies $\psi$-mixing
(see, for instance, \cite{BHK} and Ch.7,21 in \cite{Br}). Since
(\ref{2.5}) holds true with $n_0=\ell$ then we have
\begin{equation*}
P(\ell,x,\Gam)=\int_\Gam p^{(\ell)}(x,y)d\eta(y)\quad\mbox{and}\quad\mu(\Gam)=
\int_\Gam p(y)d\eta(y)
\end{equation*}
where the densities $p^{(\ell)}(x,y)$ and $p(y)$ satisfy $\gam^{-1}\geq
p^{(\ell)}(x,y)\geq\gam$ and $\gam^{-1}\geq p(y)\geq\gam$ for all
$x,y\in\cX$. Observe that (\ref{2.13+}) holds true assuming only the right
hand side of (\ref{2.5}). Indeed, since $|P(\ell,x,\Gam)-P(\ell,y,\Gam)|=
|P(\ell,x,\cX\setminus\Gam)-P(\ell,y,\cX\setminus\Gam)|$ we can assume that
$\eta(\Gam)\geq\frac 12$. Then
\[
|P(\ell,x,\Gam)-P(\ell,y,\Gam)|\leq 1-\gam\eta(\Gam)\leq 1-\frac \gam 2,
\]
and so (\ref{2.13+}) follows. Hence, the conditions of Theorem
\ref{thm2.4}(i) are satisfied and since $F_\ell=0$ $\mu^\ell$-a.e. is
excluded we conlude that the limiting variance $\sig^2$ is positive.

 Next, for each sequence $\bar x=(x^{(1)},...,x^{(\ell-1)})\in\cX^{\ell-1}$ and a
 real number $\te$ we define the (Fourier) operator
 \begin{eqnarray*}
 &\Phi_{\bar x}(\te)f(y)=E_y\exp(i\te F(x^{(1)},...,x^{(\ell-1)},\xi_\ell))
 f(\xi_\ell)\\
 &=\int P(\ell,y,dz)\exp(i\te F(x^{(1)},...,x^{(\ell-1)},z))f(z)
 \end{eqnarray*}
 and the function
 \[
 \rho_{\te}(\bar x)=\sup_{y\in\cX}\int d\eta(v)|\int p^{(\ell)}(y,z)
 p^{(\ell)}(z,v)\exp(i\te F(x^{(1)},...,x^{(\ell-1)},z))d\eta(z)|
 \]
 where $p^{(\ell)}$ is the $\ell$-step transition density of the Markov
 chain $\{\xi_n\}$. Then for any $\bar x_1=(x_1^{(1)},...,x_1^{(\ell-1)}),
 \,\bar x_2=(x_2^{(1)},...,x_2^{(\ell-1)})\in\cX^{\ell-1}$,
 \begin{eqnarray}\label{5.1}
 &\,\,\|\Phi_{\bar x_1}(\te)\Phi_{\bar x_2}(\te)\|=\sup_{f:\| f\|=1}\sup_y
 |\int d\eta(v)\exp(i\te F(x_2^{(1)},...,x_2^{(\ell-1)},v))f(v)\\
 &\times\int p^{(\ell)}(y,z)p^{(\ell)}(z,v)\exp(i\te
 F(x_1^{(1)},...,x_1^{(\ell-1)},z))d\eta(z)|\leq\rho_\te(\bar x_1).\nonumber
 \end{eqnarray}

 The following arguments aimed at proving contraction properties of the
 operators $\Phi_{\bar x}(\te)$ resemble Lemma 1.5 in \cite{Na2} though
 in our situation there are additional parameters to take care about.
 Observe that for any function $f$ and a probability measure $\nu$,
 \begin{eqnarray*}
 &|\int e^{if(x)}d\nu(x)|^2=\int\int e^{i(f(x)-f(y))}d\nu(x)d\nu(y)\\
 &=\int\int\cos(f(x)-f(y))d\nu(x)d\nu(y)=1-2\int\sin^2\frac
 {(f(x)-f(y))}2d\nu(x)d\nu(y),
 \end{eqnarray*}
 and so
 \begin{eqnarray}\label{5.1.1}
 &1-|\int e^{if(x)}d\nu(x)|\geq\frac 12(1-|\int e^{if(x)}d\nu(x)|^2)\\
 &\geq\int\sin^2\frac {(f(x)-f(y))}2d\nu(x)d\nu(y)
 =\frac 14\int\int|e^{if(x)}-e^{if(y)}|^2d\nu(x)d\nu(y)\nonumber
 \end{eqnarray}
 where we used the equality
 \[
 |e^{if(x)}-e^{if(y)}|^2=2-2\cos(f(x)-f(y))=4\sin^2\frac {(f(x)-f(y))}2.
 \]
 Now, let $\Gam_1$ and $\Gam_2$ be two subsets of the unit circle such that
 \begin{equation}\label{5.1.2}
 \inf_{\xi\in\Gam_1,\,\zeta\in\Gam_2}|\xi-\zeta|=\del>0\,\,\mbox{and}\,\,
 \min_{j=1,2}\nu\{ x:\, e^{if(x)}\in\Gam_j\}=\ve>0.
 \end{equation}
 Then for $G_j=\{ x:\, e^{if(x)}\in\Gam_j\},\, j=1,2$,
 \begin{equation}\label{5.1.3}
 \int_\cX\int_\cX |e^{if(x)}-e^{if(y)}|^2d\nu(x)d\nu(y)\geq\int_{G_1}\int_{G_2}
 |e^{if(x)}-e^{if(y)}|^2d\nu(x)d\nu(y)\geq\del^2\ve^2.
 \end{equation}

 Next, we apply the above inequalities with $\nu=\nu_{y,v}$ defined for any
  Borel $G\subset\cX$ by
 \[
 \nu(G)=\frac 1{p^{(2\ell)}(y,v)}\int_Gp^{(\ell)}(y,z)p^{(\ell)}(z,v)d\eta(z).
 \]
 Consider $g^{(\te)}_{\bar x}(y)=\exp(i\te F(\bar x,y))$, $\bar x=
 (x^{(1)},...,x^{(\ell-1)})$ as a function of $y$. Then either
 $g^{(\te)}_{\bar x}(y)$ does not depend on $y$ $\eta$-a.e. or there exist
 Borel subsets $\Gam_1=\Gam_{1,\te,\bar x}$ and $\Gam_2=\Gam_{2,\te,\bar x}$
 of the unit circle such that (\ref{5.1.2}) holds true with some
 $\del=\del_{\te,\bar x}$ and $\ve=\ve_{\te,\bar x}$. In this case by
 (\ref{2.5}) and the definition of $\nu$,
 \begin{equation}\label{5.1.4}
 \nu\{ y:\, g^{(\te)}_{\bar x}(y)\in\Gam_j\}\geq\frac {\gam^2}{p^{(2\ell)}(y,v)}
 \eta\{ y:\, g^{(\te)}_{\bar x}(y)\in\Gam_j\}\geq\frac {\gam^2\ve_{\te,\bar x}}
 {p^{(2\ell)}(y,v)}.
 \end{equation}
 This together with (\ref{5.1.1})--(\ref{5.1.3}) yields
 \begin{eqnarray}\label{5.1.5}
 &c_\te(\bar x,y,v)=p^{(2\ell)}(y,v)-|\int p^{(\ell)}(y,z)p^{(\ell)}(z,v)
 g^{(\te)}_{\bar x}(z)d\eta(z)|\\
 &=p^{(2\ell)}(y,v)(1-|\int g_{\bar x}^{(\te)}(z)
 d\nu(z)|)\geq\frac {\gam^4\del^2_{\te,\bar x}\ve^2_{\te,\bar x}}{4p^{(2\ell)}
 (y,v)}. \nonumber\end{eqnarray}
 Set $U_y=\{ v:\, p^{(2\ell)}(y,v)>2\}$. Then $\eta(U_y)\leq\frac 12
 \int_{U_y}p^{(2\ell)}(y,v)d\eta(v)\leq\frac 12$, and so
 \begin{equation}\label{5.1.6}
 \int_\cX\frac {d\eta(v)}{p^{(2\ell)}(y,v)}\geq\int_{\cX\setminus U_y}
 \frac {d\eta(v)}{p^{(2\ell)}(y,v)}\geq\frac 14.
 \end{equation}
 This together with (\ref{5.1.5}) yields that
 \begin{equation}\label{5.3}
 \rho_\te(\bar x)\leq 1-\inf_y\int c_\te(\bar x,y,v)d\eta(v)\leq
 1-\frac 1{16}\gam^4\del^2_{\te,\bar x}\ve^2_{\te,\bar x}.
 \end{equation}

 Observe that (\ref{5.3}) holds true for any $\te\ne 0$ in the non-lattice
 case if $B_{\bar x}=\emptyset$ (with $B_{\bar x}$ defined by (\ref{2.15+}))
 and for any $\te\in[-\frac \pi h,\frac \pi h]\setminus\{ 0\}$ in the lattice
  case (\ref{2.17}). Taking into account in addition that $\rho_\te(\bar x)$
  is continuous in $\te$ we conclude that for any $\del>0$ there exists
  $c_\del(\bar x)>0$ such that
  \begin{equation}\label{5.4}
  \rho_\te(\bar x)\leq 1-c_\del(\bar x)
  \end{equation}
  whenever $\del\leq |\te|\leq\frac 1\del$ and $B_{\bar x}=\emptyset$ in the
   non-lattice case and whenever $\te\in [-\frac \pi h,\frac \pi h]\setminus
   (-\del,\del)$ in the lattice case (\ref{2.17}).

   On the other hand, since $F_\ell=0$ $\mu^{\ell-1}$-a.e. is excluded then
   $g^{(\te)}_{x^{(1)},...,x^{(\ell-1)}}(y)$ cannot be constant in $y$ for
   $\mu^{\ell-1}$-almost all $(x^{(1)},...,x^{(\ell-1)})$ whenever
 $\del\leq |\te|\leq\frac 1\del$ in the non-lattice case and whenever
 $\te\in [-\frac \pi h,\frac \pi h]\setminus (-\del,\del)$ in the lattice case
 (\ref{2.17}). It follows that for each such $\te$ the inequality (\ref{5.4})
 holds true with $c_\del(\bar x)>0$ for $\bar x\in\cX^{\ell-1}$ belonging to a
 set having positive $\mu^{\ell-1}$-measure. Hence, for any $\del>0$ there
 exist $c_\del>0$ and a Borel set $G\subset\cX^{\ell-1}$ such that for all
 $\te$ from the corresponding ranges above
 \begin{equation}\label{5.5}
 \rho_\te(\bar x)\leq 1-c_\del\quad\mbox{for all}\quad\bar x\in G\,\,\,
 \mbox{and}\,\,\,\mu^{\ell-1}(G)\geq\ve>0.
 \end{equation}

    Next,
    \begin{equation}\label{5.6}
    |\vf_N(\te)|\leq E|E(\exp(i\te(S_N-S_M))|\cF_{\ell M})|
    \end{equation}
    where $M=M(N)=N-2[\frac {N-N_\ell}2]$ and, recall, $N_\ell=
    [N(1-\frac 1{2\ell})] +1$. By the Markov property we obtain that
    with probability one,
    \begin{equation}\label{5.7}
    E(\exp(i\te(S_N-S_M))|\cF_{\ell M})=\prod_{k=M+1}^N
    \Phi_{\xi_k,\xi_{2k},...,\xi_{(\ell-1)k}}(\te)\mathbf 1(\xi_{\ell M})
    \end{equation}
    where $\mathbf 1$ is the function equal 1 identically and we took into
    account that
    \begin{equation}\label{5.8}
    jM-(j-1)N\geq [\frac N2]\quad\mbox{for all}\quad j=1,2,...,\ell.
    \end{equation}
    Let $\{\xi_n^{(1)}\},\,\{\xi_n^{(2)}\},...,\{\xi_n^{(\ell-1)}\}$ be
 $\ell-1$ independent copies of the stationary Markov chain $\{\xi_n\}$ with
 the initial distribution $\mu$.
    Applying (\ref{3.1}) subsequently $\ell-1$ times with $X=(\xi_k,\,
    k\leq jN)$ and $Y=(\xi_{(j+1)n},\, n=M+1,M+2,...,N)$, $j=\ell-1,\ell-2,
    ...,1$ and using (\ref{5.8}) we obtain that
    \begin{eqnarray}\label{5.9}
    &\big\vert E|\prod_{k=M+1}^N\Phi_{\xi_k,\xi_{2k},...,\xi_{(\ell-1)k}}(\te)\mathbf 1|
    \\
    &-E|\prod_{k=M+1}^N\Phi_{\xi_k^{(1)},\xi_{2k}^{(2)},...,
    \xi_{(\ell-1)k}^{(\ell-1)}}(\te)\mathbf 1|\big\vert\leq (\ell-1)
    \psi([\frac N2]-2\ell).
    \nonumber\end{eqnarray}

     Next, introduce the Markov chain $\Xi(n)=(\xi_n^{(1)},\xi_{2n}^{(2)},...,
    \xi_{(\ell-1)n}^{(\ell-1)})$ on $\cX^{\ell-1}$ with the transition
    probability
    \[
    P_\Xi(\bar x,\bar\Gam)=\prod_{j=1}^{\ell-1}P(j,x^{(j)},\Gam_j)
    \]
    where $\bar x=(x^{(1)},...,x^{(\ell-1)})\in\cX^{(\ell-1)}$ and $\bar\Gam=
    \Gam_1\times\cdots\times\Gam_{\ell-1}$ is a Borel subset of
    $\cX^{(\ell-1)}$. Observe that $\mu^{\ell-1}=\mu\times\cdots\times\mu$
    is the stationary distribution of $\{\Xi(n)\}$ and that its transition
    probability $P_\Xi(\ell,\bar x,\cdot)$ has a transition density
    $p^{(\ell)}_\Xi$ with respect to $\eta^{\ell-1}$ satisfying
    \[
    \gam^{(\ell-1)}\leq p^{(\ell)}_\Xi(\bar x,\cdot)\leq\gam^{-(\ell-1)}
    \]
    where $\gam$ is the same as in (\ref{2.5}).

     Now, by (\ref{5.1}) and the submultiplicativity of norms of operators,
    \begin{eqnarray}\label{5.10}
    &E|\prod_{k=M+1}^N\Phi_{\xi_k^{(1)},\xi_{2k}^{(2)},...,
    \xi_{(\ell-1)k}^{(\ell-1)}}(\te)\mathbf 1|\\
    &\leq E\prod_{j=1}^{[\frac {N-N_\ell}2]}\|\prod_{k=M+2(j-1)+1}^{M+2j}
    \Phi_{\xi_k^{(1)},\xi_{2k}^{(2)},...,\xi_{(\ell-1)k}^{(\ell-1)}}(\te)\|
    \nonumber\\
    &\leq E\prod_{j=0}^{[\frac {N-N_\ell}2]}\rho_\te(\Xi(M+2j+1)).
    \nonumber\end{eqnarray}

    It follows from (\ref{2.4}), (\ref{5.1}), (\ref{5.6}), (\ref{5.7}),
    (\ref{5.9}) and (\ref{5.10}) that in order to obtain (\ref{4.16}) it
    suffices to show that for some $N_0$ and all $N\geq N_0$,
    \begin{equation}\label{5.11}
    E\prod_{n=0}^{[\frac {N-N_\ell}2]}\rho_\te(\Xi(2n+1))\leq e^{-\be N}
    \end{equation}
    for some $\be=\be_\del>0$ depending on $\del$ which determines
    corresponding domains for $\te$ in the non-lattice and the lattice cases.
    Observe that by (\ref{5.5}),
    \begin{equation}\label{5.14}
    P\{\Xi(n)\in G\}=P\{\Xi(1)\in G\}\geq\ve.
    \end{equation}

    Introduce the counting function
    \[
    V(N)=\sum^{[\frac 12(N-N_\ell)]}_{n=0}\bbI_{G}(\Xi(2n+1)),
    \]
    where $\bbI_G(x)=1$ if $x\in G$ and $=0$, otherwise, and the events
    \[
    \Gam(N)=\{ V(N)<\frac {\ve N}{9\ell^2}\}.
    \]
    Since $[\frac 12(N-N_\ell)]\geq\frac N{4\ell}-2$ it follows from
    (\ref{5.14}) and the large deviations results from \cite{Bol} together
    with \cite{DV} applied to the Markov chain $\{\Xi(2\ell n+1),\, n\geq 0\}$
    that
    \begin{equation}\label{5.15}
    P(\Gam(N))\leq\ka^{-1}e^{-\ka N}
    \end{equation}
    for some $\ka>0$ independent of $N$. Namely, by Lemma 2.5 in \cite{DV}
    the rate function on the second level of large deviations (for occupational
    measures) of the Markov chain $\Xi(2\ell n),n\geq 0$ has zero only at the
    invariant measure. According to \cite{Bol} (see also \cite{KM}) large
    deviations in the $\tau$-topology (generated by convergences on bounded
    Borel functions) have the same rate function. Thus we can apply the
    contraction principle (see \cite{DZ}) in the $\tau$-topology to conclude
    that for any bounded Borel function $g$ on $\cX^{(\ell-1)}$ the rate
    function on the first level of large deviations for sums of the form
    $\sum_{k=1}^ng(\Xi(2\ell k+1))$ has a unique zero at the integral of $g$ with
    respect to the invariant measure. Taking $g=\bbI_{G}$ we will
    arrive at (\ref{5.15}). Now by (\ref{5.4}), for all
    $\te\in[\del,\del^{-1}]$ or $\te\in[-\frac \pi{h},\frac \pi{h}]
    \setminus [-\del,\del]$ depending on the case under consideration
    \[
    \prod_{n=1}^{[\frac 12(N-N_\ell)]}\rho_\te(\Xi(n))\leq\bbI_{\Gam(N)}
    +(1-c_\del)^{\frac {N\ve}{9\ell^2}}
   \]
   and (\ref{5.11}) follows from (\ref{5.15}), completing the proof of
   Lemma \ref{lem4.1} in the setup of Theorem \ref{thm2.6}.  \qed

   Next, we derive Lemma \ref{lem4.2}. We start with the estimate (\ref{5.1})
   and the first inequality in (\ref{5.3}) but now employing the Taylor formula
   we can represent $c_\te(\bar x,y,v)$, $\bar x=(x_1^{(1)},...,x_1^{(\ell-1)})$
   there for $|\te|$ small enough in the following way (cf. Ch. 8 in
   \cite{Bre}),
  \begin{equation}\label{5.16}
  c_\te(\bar x,y,v)=\frac 12\te^2\int_\cX p^{(\ell)}(y,z)p^{(\ell)}(z,v)
  D(\bar x,y,z,v)d\eta(z)+\te^2\hat\vf_{\bar x,y,v}(\te)
  \end{equation}
   where for some constant $C>0$,
   \[
   \hat\vf_{\bar x,y,v}(\te)\leq C\vf_{\bar x}(\te)\to 0\,\,\mbox{as}\,\,
    \te\to 0\,\,\mbox{and}
    \]
 \begin{equation*}
 D(\bar x,y,z,v)=\big(F(\bar x,z)-\frac 1{p^{(2\ell)}(y,v)}\int p^{(\ell)}
 (y,z)p^{(\ell)}(z,v)F(\bar x,z)d\eta(z)\big)^2.
 \end{equation*}

 Now, either $F(\bar x,z)$ does not depend on $z$ $\eta$-a.s., i.e.
 $F_\ell(\bar x)
 =0$ $\eta$-a.s. which is excluded, or there exist Borel subsets
 $U_1=U_{1,\bar x}$ and $U_2=U_{2,\bar x}$ of the real line $\bbR$ such that
 \[
 \inf_{z\in U_1,w\in U_2}|z-w|=\del_{\bar x}>0\,\,\mbox{and}\,\,\min_{j=1,2}
 \eta(G_j(\bar x))=\ve_{\bar x}>0.
 \]
 where $G_j(\bar x)=\{ z:\, F(\bar x,z)\in U_j\},\, j=1,2$. In this case
 \begin{eqnarray}\label{5.16+}
 &\int_\cX p^{(\ell)}(y,z)p^{(\ell)}(z,v)D(\bar x,y,z,v)d\eta(z)\geq\gam^2
 \int_\cX D(\bar x,y,z,v)d\eta(z)\\
 &\geq\gam^2\inf_c\int_{G_1(\bar x)\cup G_2(\bar x)}(F(\bar x,z)-c)^2d\eta(z)
 \nonumber\\
 &\geq\gam^2\ve_{\bar x}\inf_c\inf_{a\in U_1,b\in U_2}((a-c)^2+(b-c)^2)=
 \frac 12\gam^2\ve_{\bar x}\del^2_{\bar x}.\nonumber
 \end{eqnarray}
 Now, by (\ref{5.3}), (\ref{5.16}) and (\ref{5.16+}) for $|\te|$ small enough,
 \begin{equation}\label{5.16++}
 \rho_{\te}(\bar x)\leq 1-\frac {\te^2}{4}\gam^2\ve_{\bar x}\del^2_{\bar x}+
 C\te^2\vf_{\bar x}(\te).
  \end{equation}
 Observe that by (\ref{2.1}),
 \[
 \hat D(\bar x)=\sup_y\int D(\bar x,y,z,v)d\eta(z)d\eta(v)\leq 2(1+\gam^{-3})
 \int F^2(\bar x,z)d\eta(z)<\infty\,\,\mu^{\ell-1}-\mbox{a.e.}
 \]
 Hence
 \[
 \{\bar x:\,\hat D(\bar x)\leq L\}\uparrow\cY\,\, \mbox{as}\,\, L\uparrow\infty
 \]
 where $\mu^{\ell-1}(\cY)=1$. This together with (\ref{5.16++}) yields that
 there exists $c>0$ and a Borel set $G\subset\cX^{\ell-1}$ such that
 $\mu^{\ell-1}(G)=\ve>0$ and for all $|\te|$ small enough,
 \begin{equation}\label{5.17}
 \rho_\te(\bar x)\leq e^{-c\te^2}\quad\mbox{whenever}\quad
 \bar x\in G.
 \end{equation}

 Introduce, again, the counting function
 \[
 W(N)=\sum_{n=1}^{[\frac 12(N-N_\ell)]}\bbI_G(\Xi(2n-1))
 \]
 and the events
 \[
 \Gam(N)=\{W(N)<\ve(\frac N{8\ell^2}-1)\}.
 \]
 Then, since $\rho_\te(\bar x)\leq 1$, we obtain by (\ref{5.17})
 that
 \begin{equation}\label{5.18},
 \prod_{n=1}^{[\frac 12(N-N_\ell)]}\rho_\te(\Xi(2n-1))\leq\exp(-c\ve\te^2
 (\frac N{8\ell^2}-2))+\bbI_{\Gam(N)}.
 \end{equation}
 Relying on the large deviations results from \cite{DV} and \cite{Bol} as
 explained above we conclude that
 \begin{equation}\label{5.19}
 P(\Gam(N))\leq\ka^{-1}e^{-\ka N}
 \end{equation}
 for some $\ka>0$ independent of $N$. Finally, (\ref{4.17}) follows from
 (\ref{2.4}), (\ref{5.1}), (\ref{5.6}), (\ref{5.7}), (\ref{5.9}), (\ref{5.10}),
  (\ref{5.18}) and (\ref{5.19}), completing the proof of Lemma \ref{lem4.2}
  in the setup of Theorem \ref{thm2.6}.    \qed

 \begin{remark}\label{rem5.1} Observe that weaker estimates in Lemmas
 \ref{lem4.1} and \ref{lem4.2} of the form
 $|\vf_N(\te)|\leq\frac CN$ for large $|\te|$ and $|\vf_N(\te)|\leq
 e^{-rN\te^2}+\frac CN$ for small $|\te|$, respectively, (where $C>0$ is
 a constant) also suffice for our purposes in (\ref{4.18})--(\ref{4.21}).
 Indeed, the integration in the definitions of $I_3,\, J_3,\, I_4$ and
 $J_4$ is taken over intervals with length of order $\sqrt N$, and so
 these expressions
 will still tend to zero as, first, $N\to\infty$ and then $T\to\infty$.
 But the above estimates do not require large deviations arguments since
 the latter was only needed in order to obtain exponential bounds (\ref{5.15})
 and (\ref{5.19}) for probabilities of $\Gam(N)$. On
 the other hand, in order to obtain a bound of order $\frac CN$ for these
 probabilities, which we want now, only the Chebyshev inequality is needed if
  we can show that for some constant $\tilde C>0$,
  \[
  \mbox{Var}V(N)\leq\tilde CN\,\,\mbox{and}\,\,\mbox{Var}W(N)\leq
  \tilde CN.
  \]
  It is well known (see, for instance, \cite{IL} or \cite{Br}) that the latter
  estimates will hold true if the Markov chain $\{\Xi(n)\}$ is $\alpha$-mixing
  with a summable coefficient $\al(n),\, n\geq 0$ which is much less than
  Assumption \ref{ass2.1}.
 \end{remark}

 \begin{remark}\label{rem5.2}
Most parts of the proof of Theorem \ref{thm2.6} go through assuming only the
right hand side of (\ref{2.5}) with $n_0=\ell$ replacing in the above arguments
transition densities by the Radon-Nikodim derivatives of absolutely continuos
parts of transition probabilities $P(\ell,x,\cdot)$ with respect to $\eta$
and replacing $\mu$-a.s. by $\eta$-a.s. since now these measures are not
equivalent. According to Remark \ref{rem5.1} the concluding part of the proof
 of Theorem \ref{thm2.6} needed for application of Lemmas \ref{lem4.1} and
 \ref{lem4.2} can be carried out just under $\alpha$-mixing with a summable
 coefficient which is assured when the right hand side of (\ref{2.5}) holds
 true (see \cite{Br}). Observe though that the arguments of Section \ref{sec3}
 needed to ensure positivity of the asymptotic variance $\sig^2$ as well as
 the comparison estimate (\ref{5.9}) relied on the inequality (\ref{3.1})
 which, in general, holds true only under $\psi$-mixing which is not ensured,
 in general, under just the right hand side of (\ref{2.5}). Using the technique
 from Corollary 3.6 in \cite{KV} both (\ref{5.9}) and estimates of Section
 \ref{sec3} can be caried out under substantially weaker mixing conditions
 on expense of additional moment, growth and H\" older continuity assumptions
  on the function $F$. Observe that in the case of finite Markov chains
  exponentially fast $\psi$-mixing follows from the right hand side of
  (\ref{2.5}), and so Theorem \ref{thm2.7} can be proved under this condition
  only in place of the assumptions formulated there.
\end{remark}

 \section{Finite state space case}\label{sec6}
\setcounter{equation}{0}

First, observe that our assumption that $\Pi^{k}$ for some $k$ has all
positive entries yields $\psi$-mixing of $\{\xi_n\}$ with (\ref{2.4})
satisfied (see, for instance, \cite{Br}, Ch.7). Clearly, $\Pi^n$ has all
positive entries for all $n\geq k$ and we take $m$ so that $(m-2)\ell\geq k$.
Furthermore, there exists a unique stationary distribution $\mu$ which gives
a positive mass to each state $j$ in $\cX$.
Since $F_\ell=0$ $\mu^\ell$-a.e. is excluded then Theorem \ref{thm2.4}(ii)
 guaranties positivity of
the limiting variance $\sig^2$ in (\ref{2.14}) provided $\del_\ell<1$ which
is equivalent to the first inequality in (\ref{2.20}) (see \cite{Do}).

Now, for each $\bar x=(x^{(1)},...,x^{(\ell-1)})\in\cX^{\ell-1}$ and a function
$f$ on $\cX$ we can write
\[
\Phi_{\bar x}(\te)f(a)=\sum_{b\in\cX}p^{(\ell)}_{ab}e^{i\te F(x^{(1)},...,
x^{(\ell-1)},b)}f(b).
\]
Let $\| f\|=1$ then for $\bar x_n=(x_n^{(1)},...,x_n^{(\ell-1)})\in\cX^{\ell-1},
\, n=1,...,m+2$,
\begin{eqnarray}\label{6.1}
& \|(\prod_{n=1}^{m}\Phi_{\bar x_n}(\te))f\|=\max_{a\in\cX}
|\sum_{a_1,...,a_{m}\in\cX}p^{(\ell)}_{aa_1}\\
&\times e^{i\te F(x_1^{(1)},...,x_1^{(\ell-1)},a_1)}
 p^{(\ell)}_{a_1a_2}e^{i\te F(x_2^{(1)},...,x_2^{(\ell-1)},a_2)}\nonumber\\
&\times\cdots\times p^{(\ell)}_{a_{m-1}a_{m}}e^{i\te F(x_{m}^{(1)},...,
x_{m}^{(\ell-1)},a_{m})}f(a_{m})|\leq\rho_\te(\bar x_{m-1})\nonumber
\end{eqnarray}
where for each $\bar x=(x^{(1)},...,x^{(\ell-1)})\in\cX^{\ell-1}$ we set
\begin{eqnarray*}
&\rho_\te(\bar x)=\max_a\sum_{b\in\cX}p^{((m-2)\ell)}_{ab}
\sum_{d\in\cX}|\sum_{c\in\cX}p^{(\ell)}_{bc}\\
&\times e^{i\te F(x^{(1)},...,x^{(\ell-1)},c)}p^{(\ell)}_{cd}|.
\end{eqnarray*}

Since $p^{((m-2)\ell)}_{ab}>0$ for any $a,b\in\cX$ then $\rho_\te(\bar x)=1$ if
and only if for any $b\in\cX$,
\begin{equation}\label{6.2}
\sum_{d\in\cX}|\sum_{c\in\cX}p^{(\ell)}_{bc}\\
e^{i\te F(x^{(1)},...,x^{(\ell-1)},c)}p^{(\ell)}_{cd}|=1.
\end{equation}
This will hold true if and only if $g^{(\te)}_{\bar x}(c)=
e^{i\te F(x^{(1)},...,x^{(\ell-1)},c)}$ is constant in $c$ on each set
$\cX_{b,d}=\{ c\in\cX:\, p^{(\ell)}_{bc},\, p^{(\ell)}_{cd}>0\}$.
But, in view of the assumption (\ref{2.20}) for any $c,c'\in\cX$ there exist
$b,d\in\cX$ such that $p^{(\ell)}_{bc},\, p^{(\ell)}_{cd},\, p^{(\ell)}_{bc'},
\, p^{(\ell)}_{c'd}>0$. Hence, if (\ref{6.2}) holds true then
$g^{(\te)}_{\bar x}(c)$ is constant in $c$ on the whole $\cX$.
Since we assume that $F_\ell$ is not zero identically then
   $g^{(\te)}_{x^{(1)},...,x^{(\ell-1)}}(y)$ cannot be constant in $y$ for
   $\mu^{\ell-1}$-almost all $(x^{(1)},...,x^{(\ell-1)})$ whenever
 $\del\leq |\te|\leq\frac 1\del$ in the non-lattice case and whenever
 $\te\in [-\frac \pi h,\frac \pi h]\setminus (-\del,\del)$ in the lattice case
 (\ref{2.17}). It follows that $\rho_\te(\bar x)<1$ for each such $\te$
 and for $\bar x\in\cX^{\ell-1}$ belonging to a set having positive
 $\mu^{\ell-1}$-measure. Since $\rho_\te$ is continuous in $\te$ this holds
 true uniformly in $\te$ in the corresponding compact domains above.
 Hence, for any $\del>0$ there
 exist $c_\del>0$ and a Borel set $G\subset\cX^{\ell-1}$ such that for all
 $\te$ from the corresponding ranges above
 \begin{equation}\label{6.3}
 \rho_\te(\bar x)\leq 1-c_\del\quad\mbox{for all}\quad\bar x\in G\,\,\,
 \mbox{and}\,\,\,\mu^{\ell-1}(G)\geq\ve>0.
 \end{equation}

 Next, take $M=M(N)=N-m[\frac {N-N_\ell}m]$ and proceed as in
 (\ref{5.6})--(\ref{5.9}). Introduce again the Markov chain $\Xi_n=(\xi_n^{(1)},
 \xi_{2n}^{(2)},...,\xi_{(\ell-1)n}^{(\ell-1)})$ on $\cX^{\ell-1}$ with the
 transition probabilities
    \[
    p_\Xi(\bar a,\bar b)=\prod_{j=1}^{\ell-1}p_{a^{(j)}b^{(j)}}
    \]
    where $\bar a=(a^{(1)},...,a^{(\ell-1)}),\, \bar b=(b^{(1)},...,
    b^{(\ell-1)})\in\cX^{(\ell-1)}$. If $\Pi_\Xi=(p_\Xi(\bar a,\bar b))$ is
    the corresponding transition matrix then, clearly, $\Pi_\Xi^{(m-2)\ell}$
    has all positive entries. Now similarly to (\ref{5.10}),
    \begin{eqnarray}\label{6.4}
    &E|\prod_{k=M+1}^N\Phi_{\xi_k^{(1)},\xi_{2k}^{(2)},...,
    \xi_{(\ell-1)k}^{(\ell-1)}}(\te)\mathbf 1|\\
    &\leq E\prod_{j=1}^{[\frac {N-N_\ell}m]}\|\prod_{k=M+m(j-1)+1}^{M+mj}
    \Phi_{\xi_k^{(1)},\xi_{2k}^{(2)},...,\xi_{(\ell-1)k}^{(\ell-1)}}(\te)\|
    \nonumber\\
    &\leq E\prod_{j=0}^{[\frac {N-N_\ell}m]}\rho_\te(\Xi(M+mj+1)).
    \nonumber\end{eqnarray}

    Thus, in order to obtain (\ref{4.16}) it
    suffices to show that for some $N_0$ and all $N\geq N_0$,
    \begin{equation}\label{6.5}
    E\prod_{n=0}^{[\frac {N-N_\ell}m]}\rho_\te(\Xi(mn+1))\leq e^{-\be N}
    \end{equation}
    for some $\be=\be_\del>0$ depending on $\del$ which determines
    corresponding domains for $\te$ in the non-lattice and the lattice cases.

    Introduce, again, the counting function
    \[
    V(N)=\sum^{[\frac 1m(N-N_\ell)]}_{n=0}\bbI_{G}(\Xi(mn+1)).
    \]
    where $\bbI_G(x)=1$ if $x\in G$ and $=0$, otherwise, and the events
    \[
    \Gam(N)=\{ V(N)<\frac {\ve N}{9\ell^2m}\}.
    \]
    Since $[\frac 1m(N-N_\ell)]\geq\frac N{4\ell}-2$ it follows from
    (\ref{6.3}) and the large deviations results from Section 3.1 in
    \cite{DZ} together with \cite{DV} applied to the Markov chain
    $\{\Xi(nm+1),\, n\geq 0\}$ that
    \begin{equation}\label{6.6}
    P(\Gam(N))\leq\ka^{-1}e^{-\ka N}
    \end{equation}
    for some $\ka>0$ independent of $N$. Namely, by Lemma 2.5 in \cite{DV}
    the rate function of the second level of large deviations (for occupational
    measures) of the Markov chain $\Xi(nm+1),n\geq 0$ has zero only at the
    invariant measure. The proof there depends only on the specific form of
    the rate function and it is shown in Section 3.1 in \cite{DZ} that the
    Markov chains satisfying our conditions have exactly the same form of
    the rate function. Since $\cX$ is a finite set any function is continuous
    and we can use directly the contraction principle to conclude that
    for any function $g$ on $\cX^{\ell-1}$ the rate
    function on the first level of large deviations for sums of the form
    $\sum_{k=0}^ng(\Xi(km+1))$ has a unique zero at the integral of $g$ with
    respect to the invariant measure. Now, taking $g=\bbI_{G}$ we will
    arrive at (\ref{6.6}). In fact, in our circumstances we can obtain
    (\ref{6.6}) directly from \cite{Le} and \cite{KW}. Now we conclude the
    proof of (\ref{6.5}) and of the whole Lemma \ref{lem4.1} (in the setup
    of Theorem \ref{thm2.7}) in the same way
    as in the previous section.  \qed

    Next, we derive Lemma \ref{lem4.2} in the present setup. Employing the
    Taylor formula we can write for $|\te|$ small enough (cf. Ch.8 in
    \cite{Bre}) that
    \begin{eqnarray}\label{6.7}
    &\rho_\te(\bar x)\leq 1-\frac {\te^2}2\inf_a\sum_{b\in\cX}
    p^{((m-2)\ell)}_{ab}\sum_{d\in\cX}\sum_{c\in\cX_{bd}}p^{(\ell)}_{bc}\\
&\times D(x^{(1)},...,x^{(\ell-1)}, b,c,d)p^{(\ell)}_{cd} +
\te^2\vf_{x^{(1)},...,x^{(\ell-1)}}(\te)\nonumber
\end{eqnarray}
where $\bar x=(x^{(1)},...,x^{(\ell-1)})$, $\cX_{bd}=\{ c:\, p^{(\ell)}_{bc},\,
p^{(\ell)}_{cd}>0\}$, $\vf_{x^{(1)},...,x^{(\ell-1)}}(\te)\to 0$ as $\te\to 0$
 and
\begin{eqnarray*}
&D(x^{(1)},...,x^{(\ell-1}), b,c,d)\\
&= \big (F(x^{(1)},...,x^{(\ell-1)},c)-
\frac 1{p^{(2\ell)}_{bd}}\sum_{c\in\cX_{bd}}p^{(\ell)}_{bc}p^{(\ell)}_{cd}
F(x^{(1)},...,x^{(\ell-1)},c)\big)^2
\end{eqnarray*}
while we set $\sum_{c\in\cX_{bd}}=0$ whenever $\cX_{bd}=\emptyset$. Observe
that
\begin{equation}\label{6.8}
\sum_{c\in\cX_{bd}}p^{(\ell)}_{bc}
D(x^{(1)},...,x^{(\ell-1)}, b,c,d)p^{(\ell)}_{cd}=0
\end{equation}
if and only if $F(x^{(1)},...,x^{(\ell-1)},c)$ does not depend on
$c\in\cX_{bd}$. But in view of (\ref{2.20}) for any $c,c'\in\cX$ there exist
$b,d\in\cX$ such that $c,c'\in\cX_{bd}$. It follows that (\ref{6.8}) holds
true for any $b,d\in\cX$ if and only if $F(x^{(1)},...,x^{(\ell-1)},c)$ does
not depend on $c$ on the whole $\cX$, i.e. $F_\ell(x^{(1)},...,
x^{(\ell-1)},c)=0$ for all $c\in\cX$. Since the latter equality cannot hold
true identically by our assumptions we conclude that there exists a nonempty
 $G\subset\cX^{\ell-1}$ such that for some $r>0$ and all $\bar x\in G$,
 \begin{equation}\label{6.9}
 \rho_\te(\bar x)\leq 1-r\te^2\leq e^{-r\te^2}
 \end{equation}
 provided $|\te|$ is small enough. In present circumstances
 $\mu^{\ell-1}(G)>0$ as for any nonempty subset of $\cX^{\ell-1}$ and
 we complete the proof of Lemma \ref{lem4.2} (in the present setup) proceeding
  in the same way as in
 the proof of Lemma \ref{lem4.1} above and at the end of the previous section
 by introducing a counting function of arrivals of $\Xi(mn+1)$ to $G$,
 using (\ref{6.4}) and (\ref{6.9}), and relying on the same large deviations
 argument. \qed

\end{document}